\documentclass[a4paper,10pt]{article}

\usepackage{amssymb}

\usepackage{amsmath}

\usepackage{tikz}
\def\NZQ{\Bbb}               
\def\NN{{\NZQ N}}
\def\KK{{\NZQ K}}

\def\Bc{{\mathcal{B}}}
\def\Fc{{\mathcal{F}}}

\def\Mc{{\mathcal{M}}}
\def\Vc{{\mathcal{V}}}

\def\ab{{\bold a}}
\def\kb{{\bold k}}
\def\xb{{\bold x}}

\def\mm{{\mathfrak m}}

\def\mdeg{{\rm mdeg}}
\def\lcm{{\rm lcm}}
\def\low{{\rm low}} 
\def\LFD{{\rm LFD}}
\def\GFD{{\rm GFD}}
\def\MVT{{\rm MVT}}

\def\cocoa{\hbox{\rm C\kern-.13em o\kern-.07em C\kern-.13em o\kern-.15em A}}
\def\cocoalib{\hbox{\rm C\kern-.13em o\kern-.07em C\kern-.13em o\kern-.15em A\kern-.13em
	L\kern-.13em i\kern-.13em b}}

\newtheorem{Theorem}{Theorem}[section]

\newtheorem{definition}[Theorem]{Definition}
\newtheorem{example}[Theorem]{Example}
\newtheorem{remark}[Theorem]{Example}
\newtheorem{proposition}[Theorem]{Proposition}
\newtheorem{lemma}[Theorem]{Lemma}
\newtheorem{corollary}[Theorem]{Corollary}

\title{(n-1)-st Koszul homology and the structure of monomial ideals}

\author{Anna M. Bigatti \thanks{Universitá degli Studi Genova (Italy)
              {\tt bigatti@dima.unige.it}}
		 \and Eduardo S\'aenz-de-Cabez\'on
		\thanks{ Universidad de La Rioja (Spain)
              {\tt eduardo.saenz-de-cabezon@unirioja.es}}}

\begin{document}


\maketitle

\begin{abstract}
Koszul homology of monomial ideals provides a description of the structure of such ideals, not only from a homological point of view (free resolutions, Betti numbers, Hilbert series) but also from an algebraic viewpoint. In this paper we show that, in particular, the homology at degree $(n-1)$, with $n$ the number of indeterminates of the ring, plays an important role for this algebraic description in terms of Stanley and irreducible decompositions.
\end{abstract}

\section{Introduction}
\label{intro}
Koszul homology of a module over a polynomial ring is a very relevant invariant that describes the structure of such modules \cite{B64,E95,K50}, we can read from it the main homological informations about the module, namely its Betti numbers, minimal free resolutions, Hilbert series, etc. In the case of monomial ideals, their combinatorial nature is reflected in their Koszul homology \cite{MS04}, we can then speak of \emph{combinatorial Koszul homology} \cite{S08}.

In this context, combinatorial Koszul homology can also be used to describe the \emph{algebraic} structure of monomial ideals, in terms of irreducible and Stanley decompositions. The $(n-1)$-st Koszul homology modules play here a very relevant role due to the correspondence between the multidegrees in which the corresponding $(n-1)$-st Koszul homology is not null and the maximal standard monomials with respect to the ideal. This correspondence allows us not only to describe the above named decompositions of a monomial ideal but also to compute them modifying existing algorithms used to compute Koszul homology of monomial ideals. This kind of algorithms show good performance when compared to algorithms specialized in irreducible decompositions.

The paper is organized as follows: Section \ref{sec:homology} gives the definition and basic properties of the Koszul homology of monomial ideals. We give also here some notations that will be used through the text. Sections \ref{sec:irreducible} and \ref{sec:stanley} contain the main results and procedures to obtain algebraic descriptions of monomial ideals from their Koszul homology. Finally, in Section \ref{sec:algorithms} we show the algorithms that we have implemented to obtain irreducible decompositions using their relation to Koszul homology, and compare with other approaches.

\section{Koszul homology of monomial ideals}
\label{sec:homology}
Let $R=\kb[x_1,\dots,x_n]$ be the polynomial ring in $n$ indeterminates over a field $\kb$ of caracteristic $0$. Let $\Vc$ be an $n$ dimensional $\kb$-vector space, and $S\Vc$, $\wedge\Vc$ the symmetric and exterior algebras over $\Vc$ respectively. Let us take a basis of $\Vc$ and denote it $\{x_1,\dots,x_n\}$ so that we can identify $S\Vc$ and $R$.

Consider the following complex:
$$\KK: 0\rightarrow R\otimes_\kb\wedge^n\Vc\stackrel{\partial}{\rightarrow}R\otimes_\kb\wedge^{n-1}\Vc\stackrel{\partial}{\rightarrow}\cdots R\otimes_\kb\wedge^{1}\Vc\stackrel{\partial}{\rightarrow}R\otimes_\kb\wedge^0\Vc\rightarrow 0$$

where the differential $\partial$ is defined by the following rule

$$\partial(x_1^{\mu_1}\dots x_n^{\mu_n}\otimes x_{j_1}\wedge\cdots\wedge x_{j_i})=\sum_{k=1}^i (-1)^{k+1} x_{j_k}\cdot x_1^{\mu_1}\dots x_n^{\mu_n}\otimes x_{j_1}\wedge\cdots\wedge \hat{x_{j_k}}\wedge\cdots \wedge x_{j_i}$$

This differential verifies $\partial^2=0$ and makes $\KK$ a complex, which is called the {\bf Koszul complex}. This complex is exact and it is therefore a minimal free resolution of the base field $\kb$.

\begin{definition}
Let $\Mc$ be an $R$-module. The complex $\KK(\Mc):=\Mc\otimes_R\KK$ is called the \emph{Koszul complex} of $\Mc$. $\KK(\Mc)$ is not exact in general, and its homology is called the \emph{Koszul homology} of $\Mc$; it will be denoted $H_*(\KK(\Mc))$ or $H_*(\Mc)$ for short.
\end{definition}

If the module $\Mc$ is (multi)-graded then $\KK(\Mc)$ and $H_*(\Mc)$ are also (multi)-graded, since $\partial$ preserves multidegree. Monomial ideals are a particular case of multigraded modules with the natural multigrading $\mdeg(x_i)=(0,\dots,\buildrel{i}\over{1},\dots,0)$. This multigrading induces the following multigrading in $\KK(\Mc)$: $\mdeg(x^\mu\otimes x_{j_1}\wedge\cdots\wedge x_{j_i}))=\mu+(0,\dots,\buildrel{j_1}\over{1},\dots,\buildrel{j_k}\over{1},\dots,0)$. We denote $H_{i,\mu}(\KK(I))$ the multidegree $\mu$ component of $H_i(\KK(I))$.

\begin{example}\label{ex:algebraic_n-1}
Let $I\subseteq\kb[x_1\dots,x_n]$ be a monomial ideal. Let $c\in I\otimes \wedge^{n-1}\Vc$ such that $\mdeg(c)=\mu$. Then $c$ is of the form:
$$c=a_1\otimes x_2\wedge\dots\wedge x_n+a_2\otimes x_1\wedge x_3\wedge\dots\wedge x_n+a_n\otimes x_1\wedge\dots\wedge x_{n-1}$$
where the $a_i$ are monomials. Since $\mdeg(c)=\mu$ then $\mdeg(a_i)=\mdeg(\frac{x_i\cdot x^\mu}{x_1\cdots x_n})$. Now, $c$ is a cycle iff $\partial(c)=0$. On the other hand, $c$ is a boundary if there exists $c'\in I\otimes \wedge^{n}\Vc$ such that $\partial c'=c$. If $c'=b\otimes x_1\wedge\cdots\wedge x_n$ then 
$$\partial(c')=x_1\cdot b\otimes x_2\wedge\dots\wedge x_n-x_2\cdot b\otimes x_1\wedge x_3\wedge\dots\wedge x_n+(-1)^{n+1}x_n\cdot b\otimes x_1\wedge\dots\wedge x_{n-1}.$$
Since $\partial$ preserves multidegree, we have that $\mdeg(c')=\mdeg(c)$, hence
$$\mdeg(b)=\mdeg(\frac{x^\mu}{x_1\cdots x_n}).$$
Taking now a monomial $x^\mu\in I$, if $\frac{x_i\cdots x^\mu}{x_1\cdots x_n}\in I$ we take $a_i=(-1)^{i+1}\frac{x_i\cdots x^\mu}{x_1\cdots x_n}$, then $c=\sum_{i=1}^n a_i\otimes x_1\wedge\cdots\wedge\hat{x_i}\wedge\cdots\wedge x_n$ is a cycle in $I\otimes\wedge^{n-1}\Vc$ which is a boundary iff $\frac{x^\mu}{x_1\cdots x_n}\in I$. If $\frac{x^\mu}{x_1\cdots x_n}\notin I$ then $c$ is a generator of $H_{n-1,\mu}(\KK(I))$. It is not hard to see that it is in fact the only generator. Hence
$$H_{n-1,\mu}(\KK(I))\simeq\left\{ \begin{array}{lr}\kb&\mbox{ if } \frac{x_i\cdot x^\mu}{x_1\cdots x_n}\in I \,\forall i \mbox{ and } \frac{x^\mu}{x_1\cdots x_n}\notin I\\0&\mbox{ in any other case}\end{array}\right.$$

Later in example \ref{ex:simplicial_n-1} we give a simplicial version of this fact.
\end{example}

\begin{remark}
A basic equality between the $\kb$-vector space dimension of the Koszul homology modules and the Betti numbers of a monomial ideal is based on two equivalent ways of computing certain $Tor$ modules. This equality will be underlying all the work in the following pages:
\begin{equation}\label{tor}
 dim_{\kb}(\KK_{i,\mu}(I))=dim_{\kb}(Tor_{i,\mu}(I,\kb))=\beta_{i,\mu}(I)
\end{equation}
\end{remark}

To finish this section we give some notations that will be used through the paper and we also introduce two simplicial complexes that can be associated to a monomial ideal $I$ and a multidegree $\mu$ the simplicial homology of which is equivalent to the Koszul homology of $I$ at multidegree $\mu$. These complexes are a prominent example of techniques at the interplay of combinatorics and algebra that are ubiquitous when dealing with monomial ideals.

\paragraph*{Notations: }
\begin{itemize}
\item Given a multidegree $\mu\in \NN^n$ we denote by $supp(\mu)$ the nonzero indices of $\mu$ and call it the \emph{support} of $\mu$. The support of the monomial $x^\mu$ is the set $\{x_i\vert i\in supp(\mu)\}$. We say that a multidegree $\mu$ (correspondingly a monomial $x^\mu$) has full support if $supp(\mu)=\{1,\dots,n\}$.
\item For every monomial ideal $I$, we denote by $\Bc_{n-1}(I)$ the set of multidegrees $\mu$ such that $H_{n-1,\mu}(I)\neq 0$.
\end{itemize}

\begin{definition}
Let $\mu\in\NN^n$, we say that its \emph{lowered} multidegree $\low(\mu)$ is $\mu'$ where $\mu'_i=max(\mu_i-1,0)$ i.e. we
substract one from every index in the support of $\mu$.
\end{definition}

\begin{definition}
Let $I\subseteq \kb[x_1,\dots,x_n]$ a monomial ideal, and $\mu\in\NN^n$. We define the (upper) \emph{Koszul simplicial complex} associated to $I$ at $\mu$, denoted $\Delta^I_\mu$ as
$$\Delta^I_\mu=\{\mbox{squarefree vectors }\tau\vert x^{\mu-\tau}\in I\}$$
Dually, we define the \emph{lower Koszul simplicial complex} associated to $I$ at $\mu$, denoted $\Delta^\mu_I$ as
$$\Delta^\mu_I=\{\mbox{squarefree vectors }\tau\vert x^{\low(\mu)+\tau}\notin I\}$$
\end{definition}

The relation between the Koszul homology of $I$ and the simplicial homology of these complexes is the following:

\begin{proposition}
\begin{equation}
 H_{i,\mu}(\KK(I))\simeq \tilde{H}_{i-1}(\Delta^I_{\mu})\simeq\tilde{H}^{\vert supp(\mu)\vert-i-2}(\Delta^\mu_I)
\end{equation}
in this context $\tilde{H}^i$ simply means $\tilde{H}_{-i}$.
\end{proposition}
\paragraph{Proof: } For the first equivalence, we use the isomprphism
\begin{equation}\label{simplicial_isomorphism}
\tau \mapsto\frac{x^\mu}{x^\tau}\otimes x^{\tau_1}\wedge\cdots\wedge x^{\tau_n}\,\forall \tau=(\tau_1,\dots,\tau_n).
\end{equation}
The second equivalence is just the Alexander duality between $\Delta^I_\mu$ and $\Delta^\mu_I$ (see \cite{MS04}).

%
%

\begin{example}\label{ex:simplicial_n-1}
Let $I\subseteq\kb[x_1,\dots,x_n]$ a monomial ideal. We have that $H_{n-1,\mu}(\KK(I))\simeq \tilde{H}_{n-2}(\Delta^I_{\mu})$. Since $\Delta^I_\mu$ is a subcomplex of the standard $n$-simplex $\Delta_n$, we have that the only possibility of $\Delta^I_\mu$ has homology at degree $n-2$ is that all $(n-1)$-faces of $\Delta_n$ are in $\Delta^I_\mu$ and the only $n$-face of $\Delta_n$ is not in $\Delta^I_\mu$, in which case $H_{n-2}(\Delta^I_\mu)\simeq\kb$. Observe that due to (\ref{simplicial_isomorphism}) the $(n-1)$- faces correspond to the monomials $\frac{x_i\cdot x^\mu}{x_1\cdots x_n}$ and the $n$-face corresponds to $\frac{x^\mu}{x_1\cdots x_n}$. Therefore: 
$$H_{n-1,\mu}(\KK(I))=H_{n-2}(\Delta^I_\mu)\simeq\left\{ \begin{array}{lr}\kb&\qquad\mbox{if } \frac{x_i\cdot x^\mu}{x_1\cdots x_n} \in I\,\forall i \mbox{ and } \frac{x^\mu}{x_1\cdots x_n}\notin I\\0&\mbox{ in any other case}\end{array}\right.$$
\end{example}

From the example \ref{ex:simplicial_n-1} we can see that the geometric realization of the multidegrees $\mu\in\Bc_{n-1}(I)$ can be called \emph{maximal corners}. We shall give the same name to the corresponding multidegree:

\begin{definition}
We say that a monomial $x^\mu$ is a \emph{closed corner} if $x_i\cdot \mu'\in I\, \forall x_i\in supp(\mu)$. We say that a monomial is a \emph{maximal corner} if it is a closed corner with maximal support. 
\end{definition}

\begin{figure}
\begin{center}
\begin{tikzpicture}[scale=1]
\draw [very thin, ->](0,0)--(-0.7,-0.7);
\draw [very thin, ->](4.5,1.5)--+(1,0);
\draw [very thin, ->](1.5,4.5)--+(0,1);
\foreach \position in {(0,0),(1.5,0.5),(2.5,0.5),(1,2),(1,3),(2,2),(2,3),(3,2),(3,3)} 
\filldraw[fill=gray!65,draw=black]\position rectangle +(1,1);
\foreach \position in {(1,0),(3.5,0.5),(4,1),(4,2),(4,3)}
\filldraw[fill=gray!40!white,draw=black] \position-- ++(0.5,0.5)-- ++(0,1)-- ++(-0.5,-0.5)-- cycle;
\foreach \position in {(0,1),(0.5,1.5),(1.5,1.5),(2.5,1.5),(1,4),(2,4),(3,4)}
\filldraw[fill=gray!15!white,draw=black] \position-- ++(1,0)-- ++(0.5,0.5)-- ++(-1,0)-- cycle;
\fill (0,0) circle(2pt) node[left]{$x^3$};
\fill (1.5,0.5) circle(2pt) node[below right]{$x^2y$};
\fill (1,2) circle(2pt) node[left]{$xz$};
\fill (4.5,1.5) circle(2pt) node[above right]{$y^3$};
\fill (1.5,4.5) circle(2pt) node[above left]{$z^3$};

\fill[color=white, draw=black](4,4) circle(2pt) node[below left]{$xy^3z^3$};
\fill[color=white,draw=black](3.5,1.5) circle(2pt) node[below left]{$x^2y^3z$};
\fill[color=white,draw=black](1,1) circle(2pt) node[below left]{$x^3yz$};

\fill[color=black!50,draw=black](4.5,4.5) circle(2pt) node[above right]{$y^3z^3$};
\fill[color=black!50,draw=black](3.5,0.5) circle(2pt) node[below right]{$x^2y^3$};
\fill[color=black!50,draw=black](0,1) circle(2pt) node[above left]{$x^3z$};
\fill[color=black!50,draw=black](1,0) circle(2pt) node[below right]{$x^3y$};
\fill[color=black!50,draw=black](1,4) circle(2pt) node[above left]{$xz^3$};

\end{tikzpicture}
\caption{Closed and maximal corners of the ideal $I=\langle x^3,x^2y,xz,y^3,z^3\rangle$}\label{fig:closed-maximal}
\end{center}
\end{figure}
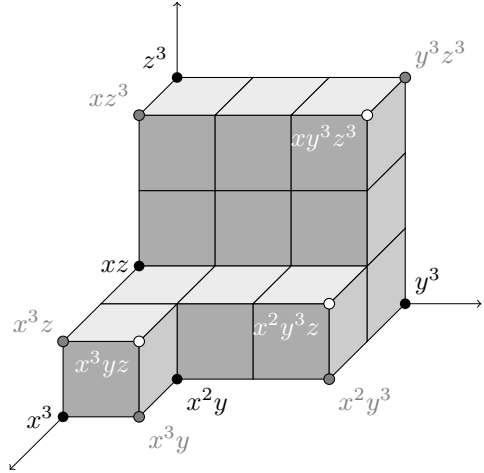

It is clear that $x^\mu$ is a maximal corner if and only if $\mu\in\Bc_{n-1}(I)$, for general closed corners only the direct statement is verified i.e. $x^\mu$ is a closed corner $\Rightarrow H_{\vert supp(x^\mu)\vert -1,\mu}(\KK(I))\neq 0$. Maximal corners correspond to maximal standard monomials, more precisely $x^\mu$ is a closed corner iff $x^{\mu'}$ is a maximal standard monomial, i.e. a monomial that is not in $I$ and is maximal with respect to divisibility among the monomials not in $I$. Figure \ref{fig:closed-maximal} shows an example of maximal and closed corners.

We finish the section with some definitions that will be needed in the rest of the paper.

\begin{definition}
\begin{itemize}
\item Let $I$ be a monomial ideal, its \emph{artinian closure} $\hat{I}$ is the smallest artinian ideal such that the set of minimal generators of $I$ is included in the set of minimal generators of $\hat{I}$. 
\item If $I=\langle m_1,\dots,m_r\rangle$ then $\lambda(I):=\lcm(m_1,\dots,m_r)$ denotes the least common multiple of all the minimal generators of $I$.
\end{itemize}

\end{definition}
Let $\lambda= \lambda(I)$, observe that  $\hat{I}=I+\langle x_1^{\lambda_1+1},\dots,x_n^{\lambda_n+1}\rangle$. For any artinian monomial ideal we have $I=\hat{I}$.

\section{Koszul homology and Irreducible decompositions}\label{sec:irreducible}

We say that a monomial ideal is \emph{irreducible} if it is generated by powers of the indeterminates (cf. def. 5.6.16 and prop. 6.2.11 in \cite{KR05} for a proof of this caracterization of irreducible monomial ideals). Given $a\in\NN^n$ we indicate with $\mm^\ab$ the irreducible monomial ideal $\langle x_i^{a_i} \mid a_i\ne 0\rangle$. The ideal $\mm=\langle x_1,\dots,x_n\rangle$ is called the \emph{irrelevant} ideal. Every monomial ideal has a unique irredundant irreducible decomposition \cite{V01}. The irredundant irreducible decomposition of any monomial ideal can be obtained in terms of its Koszul homology. First we give some preparatory results.

\begin{lemma}\label{lem:union}
Let $x^\mu$ a monomial and let $C$ be the complement of $supp(x^\mu)$ in the set of all indeterminates. Then, for any monomial $x^\nu$ we have
$$x^\nu\notin \mm^\mu \iff x^\nu \in \bigcup_{x^\rho\vert x^{\low(\mu)}}x^\rho\cdot\kb[C]$$
\end{lemma}

\paragraph{Proof: } Let $\mu'=\low(\mu)$. For the direct statement assume $x^\nu\not\in\mm^\mu$. Then $\nu_i\lneq\mu_i$ for all $\mu_i\neq 0$, i.e. for all $i\notin C$, then clearly $x^\nu \in \bigcup_{x^\rho\vert x^{\mu'}}x^\rho\cdot\kb[C]$.

For the other direction let $x^\nu\in\mm^\mu$, then $\nu_i\geq\mu_i$ for some $i\notin C$, thus $x^\nu \notin \bigcup_{x^\rho\vert x^{\mu'}}x^\rho\cdot\kb[C]$.

\begin{corollary}\label{col:fullsupport}
Let $x^\mu$ and $x^\nu$ be two different monomials such that $x^\mu$ has full support. Then
$$x^\nu\notin\mm^\mu\iff x^\nu\vert x^{\low(\mu)}$$
\end{corollary}

\paragraph{Proof: }In this case $C=\emptyset$.

Now we are ready to give an expression of the irredundant irreducible decomposition of $I$. Let us distinguish two cases:

\begin{proposition}\label{prop:irred_1}
Let $I$ be an \emph{artinian} monomial ideal and let $\Bc_{n-1}(I)$ the set of multidegrees in which $I$ has nonzero $(n-1)$-st Koszul homology (i.e. the maximal corners of $I$ ). The irredundant irreducible decomposition of $I$ is
$$I=\bigcap_{\mu\in\Bc_{n-1}(I)}\mm^\mu$$
\end{proposition}

\paragraph{Proof: } Since $I$ is artinian $R/I$ is a finite $\kb$-vector space, i.e. the set of standard monomials is finite. Every standard monomial divides a maximal standard monomial by definition, and the set of maximal standard monomials is given by the monomials $x^{\low(\mu)}$ such that $\mu\in\Bc_{n-1}(I)$.

Applying now Corollary \ref{col:fullsupport}, we have that $x^\nu\vert x^{\low(\mu)}\iff x^\nu\notin\mm^\mu$, therefore $x^\nu\notin I\iff x^\nu\in R/I\iff x^\nu\notin \bigcup_{\mu\in\Bc_{n-1}(I)}\mm^\mu$ i.e. $x^\nu\in I\iff x^\nu\in \bigcap_{\mu\in\Bc_{n-1}(I)}\mm^\mu$ and we have an irreducible decomposition of $I$. Since the elements in $\Bc_{n-1}(I)$ do not divide each other, the decomposition is irredundant.

\begin{proposition}
Let $I$ be a monomial ideal, $\hat{I}$ its artinian closure and $\lambda=\lambda(I)$. The irredundant irreducible decomposition of $I$ is given by

 $$I = \bigcap_{\mu\in\Bc_{n-1}(\hat{I})}\mm^{\tilde{\mu}}  \qquad{\rm where}\quad
  \tilde{\mu}_i=
  \begin{cases}
    \mu_i & {\rm if\ } \mu_i \le \lambda_{i}\\
    0     & {\rm otherwise}
  \end{cases}
 $$
\end{proposition}

\paragraph{Proof: }

Let $x^\mu$ be a maximal corner of $\hat{I}$ and $C_\mu$ the set of indeterminates such that $\mu_i\geq\lambda_i +1$.

We have from Lemma \ref{lem:union} that 
$$x^\nu\notin\bigcap_{\mu\in\Bc_{n-1}(\hat{I})}\mm^{\tilde{\mu}}\iff x^\nu\in\bigcup_{\stackrel{x^\rho\vert x^{\low(\tilde{\mu})}}{\mu\in\Bc_{n-1}(\hat{I})}}x^\rho\cdot\kb[C_\mu]$$
Let us call $J$ the set on the right hand side of this equivalence.

Consider now a monomial $x^\nu$ in $J$ such that $x^\nu$ is not in $\hat{I}$, then clearly $x^\nu\notin I$. Observe that in particular $x^{\low(\mu)}\notin I$.

Take then $x^\nu\in J$ such that $x^\nu\in \hat{I}$. We have that $x^\nu$ is of the form $x^\rho\cdot\kb[C_\mu]$ with $x^\rho\vert x^{\low(\tilde{\mu})}$ for some $\mu\in\Bc_{n-1}(\hat{I})$. Since $x^\nu\in\hat{I}$, $C_{\mu}\neq\emptyset$, thus $\exists \nu_i\leq\lambda_i+1$. Assume $x^\nu\in I$, then there is some $m\in\min(I)$ such that $m_i\leq \nu_i$ for all $i$. In particular $m_i\leq \nu_i$ for all $i\notin C_\mu \Rightarrow m\vert x^{\low(\mu)}\Rightarrow x^{\low(\mu)}\in I$ which is a contradiction.

Therefore, $J\subseteq R/I$ (the inclusion $R/I\subseteq J$ is trivial), then $J=R/I$, or equivalently
$$I=\bigcap_{\mu\in\Bc_{n-1}(\hat{I})}\mm^{\tilde{\mu}}$$

\begin{remark}
Given any monomial ideal $I$, a procedure is given in \cite{S08} to obtain the irredundant irreducible decomposition of $I$ directly from its Koszul homology. In this case one needs to know the multidegrees of the generators of every homology module (or the full minimal free resolution), not only the $(n-1)$-st module.
\end{remark}

\begin{example}

Consider the ideal $I=\langle x^3y^5z,y^5z^4,y^3z^5,xyz^5,x^2z^5,x^4z^3,x^4y^2z^2,x^4y^4z\rangle$. The {\lcm} of its minimal generators is $\lambda(I)=(x^4y^5z^5)$. The maximal corners of $\hat{I}$ are $x^4y^5z^5, x^5y^2z^3,x^5y^4z^2,x^3y^6z^4,x^2yz^6,xy^3z^6$ and $x^5y^6z$.
Thus, we obtain the irredundant irreducible decomposition:
$$I=\langle x^4,y^5,z^5\rangle\cap\langle y^2,z^3\rangle\cap\langle y^4,z^2\rangle\cap\langle x^3,z^4\rangle\cap\langle x^2,y\rangle\cap\langle x,y^3\rangle\cap \langle z\rangle$$
\end{example}

\section{Koszul homology and Stanley decompositions}\label{sec:stanley}
Given an ideal $I\subseteq R$, the quotient ring $R/I$ is a $\kb$-algebra. If $I$ is a monomial ring, the monomials in $R/I$ can be seen as the complement of the monomials in $I$ in the set of all monomials, and they form a  $\kb$-vector space basis of $R/I$. We call combinatorial decomposition of $R/I$ to a representation of it as a finite sum of $\kb$-vector spaces of the form $x^\mu\kb[\xb_\mu]$ with $\xb_\mu\subseteq[x_1,\dots,x_n]$. If the sum is direct we say it is a \emph{Stanley decomposition}:
\begin{equation}\label{eq:stanley}
 R/I=\bigoplus_{\mu\in\Fc}x^\mu\kb[\xb_\mu]
\end{equation}
where $\Fc$ is a finite subset of $\NN^n$. Note that Stanley decompositions are not unique. Every ideal $I\subseteq R$ has a Stanley decomposition. This type of decompositions provide information on some relevant invariants of $R/I$:
\begin{proposition}[S78,SW91]
Let \emph{(\ref{eq:stanley})} be a Stanley decomposition of $R/I$, and let $d$ be the maximum of the numbers $\vert \xb_\mu\vert$, $\mu\in\Fc$. Then
\begin{enumerate}
 \item $d$ is the Krull dimension of $R/I$.
 \item The Hilbert series of $R/I$ is given by
	$$H(R/I;t)=\sum_{\mu\in\Fc}\frac{t^{\vert\mu\vert}}{(1-t)^{\vert\xb_\mu\vert}}$$
\end{enumerate}
\end{proposition}

To explain how to obtain Stanley decompositions of monomial ideals from their Koszul homology, we need some auxiliary concepts:

\begin{definition}
We say that a set of indeterminates $\{x_{j_1},\dots,x_{j_k}\}$ is a \emph{cone of locally free directions} of the monomial $x^\mu$ with respect to $I$ if $\tau=\{j_1,\dots,j_k\}\in\Delta^\mu_I$. The set of cones of locally free directions of $x^\mu$ will be denoted $\LFD(x^\mu)$, and is given by the facets of $\Delta^\mu_I$.

We say that a set of indeterminates $\{x_{j_1},\dots,x_{j_k}\}$ is a \emph{cone of globally free directions} of the monomial $x^\mu\notin I$ with respect to $I$ if $x^\mu\cdot x^\sigma\notin I$ for all monomials $x^\sigma\in\kb[x_{j_1},\dots,x_{j_k}]$. The set of cones of globally free directions of $x^\mu$ will be denoted $\GFD(x^\mu)$. For $x^\mu\in I$ we state $\GFD(x^\mu):=\GFD(x^{\low(\mu)})$.
\end{definition}

\begin{example}
Consider the ideal $I=\langle x^3,x^2y,xz,y^3,z^3\rangle$.

We have that $\Delta^{xyz}_I=\{\{x\},\{y\},\{z\},\{x,y\},\{y,z\}\}$ thus,  $\LFD(xyz)=\{[x,y],[y,z]\}$. However, $\GFD(xyz)=\emptyset$ since $\low(xyz)=1$ and $1\cdot x^3\in I$, $1\cdot y^3\in I$ and $1\cdot z^3\in I$.

Note that if we consider instead $I=\langle x^2y,xz,y^3,z^3\rangle$ then $\LFD(xyz)$ remains unchanged, but $\GFD(xyz)=\{[x]\}$ since $1\cdot x^p\notin I$ for all $p\geq 0$.
\end{example}

We can now proceed to describe how to obtain a Stanley decomposition of a monomial ideal from its Koszul homology. We treat separately the artinian and non-artinian cases.

\begin{proposition}
Let $I$ be an artinian monomial ideal. A Stanley decomposition of $R/I$ is given by
$$R/I\simeq \bigoplus_{\substack{x^\nu\vert x^{\low(\mu)}\\ x^\mu \in \Bc_{n-1}(I)}} \kb\cdot x^\nu$$
\end{proposition}

\paragraph{Proof: }
We just need to remind the first part of the proof of Proposition \ref{prop:irred_1}: Since $I$ is artinian $R/I$ is a finite $\kb$-vector space, i.e. the set of standard monomials is finite. Every standard monomial divides a maximal standard monomial by definition, and the set of maximal standard monomials is given by the monomials $x^{\low(\mu)}$ such that $\mu\in\Bc_{n-1}(I)$. 

\begin{proposition}
Let $I$ a monomial ideal and let $\hat{I}$ its artinian closure. There is a procedure to obtain a Stanley decomposition of $I$ form the $(n-1)$-st Koszul homology of $\hat{I}$.
\end{proposition}

\paragraph{Proof: } The proof consists on a description of such procedure.
\begin{enumerate}
\item First we compute the Stanley decomposition of $R/\hat{I}$ using the previous proposition. Since $R/\hat{I}\subseteq R/I$, this is part of the Stanley decomposition of $R/I$, we call it the \emph{inner part} of the decomposition.
\item Take now $x^\mu\in \Bc_{n-1}(\hat{I})$ such that $\mu_i\geq \lambda_i+1$ for some $i\in\{1,\dots,n\}$, where $\lambda=\lambda(I)$. We call \emph{points of the skeleton} (of the decomposition of $R/I$) to the monomials $\frac{x^\tau\cdot x^\mu}{x_1\cdots x_n}$ such that $\tau \in \Delta^I_\nu$ and also to the monomials that are divisors of $\frac{x^\tau\cdot x^\mu}{x_1\cdots x_n}$ in the nonfree directions of $\{x_1,\dots,x_n\}\setminus\tau$
\item To obtain the Stanley decomposition we add the cones of the points of the skeleton in all their free directions. 
\end{enumerate}

Observe that if a point is a {\it nonfree divisor} of several points in the skeleton, then all these point have the same free directions. To prove this, assume $x^\rho$ divides both $x^\mu$ and $x^\nu$ and let $GFD(x^\mu)$ be the set of free directions of $x^\mu$ and $GFD(x^\nu)$ the set of free directions of $x^\nu$. We know that $GFD(x^\mu)\cap supp(\mu-\rho)=\emptyset$ and $GFD(x^\nu)\cap supp(\nu-\rho)=\emptyset$. Assume now that there exists $i\in GFD(x^\mu)-GFD(x^\nu)$, then we must have $\rho_i=\mu_i=\lambda_i+1$ and $\nu_i<\lambda_i+1$ and we have a contradiction, because in that case $x^\rho$ would not be a divisor of $x^\nu$.

Now, let $x^\mu\in R/I$, then we can have two cases, first, if $x^\mu$ divides $x^{\lambda+1}$ then it is in the inner part of $\hat{I}$, and from the first considerations above we know that we have collected all the inner part of the Stanley decomposition of $R/I$. Second, if $x^\mu$ does not divide $x^{\lambda+1}$ then exists $i$ such that $\mu_i>\lambda_i+1$, for each such $i$ do $\bar{\mu}=\mu-(\mu_i-\lambda_{i}+1))_i$, and then we are back in the first case at a point in which all $i$ are free directions, and we know it is in some of the considered cones.

\section{Algorithms}\label{sec:algorithms}
We present in this section an algorithm for computing irreducible decompositions of monomial ideals based on the correspondence between the irreducible components of the monomial ideal $I$ and the multidegrees in $\Bc(I)$. The algorithm is a modification of the one used to compute Koszul homology using Mayer-Vietoris trees presented in \cite{S08}.

The components of the irredundant irreducible decomposition of a monomial ideal $I$ also correspond to the minimal generators of its Alexander dual \cite{MS04}, to the facets of its Scarf complex \cite{M04}, and to the maximal standard monomials associated to $I$ \cite{R08}. These correspondences have originated different algorithms for the computation of irredundant irreducible decompositions. The most relevant ones are those proposed by Milowski \cite{M04}, based on Alexander duals and Scarf complexes; and Roune, based on enumeration of maximal standard monomials \cite{R08}. This last approach is the closest to our approach by $(n-1)$-st Koszul homology of $I$.

\subsection{Mayer-Vietoris trees}
Given a monomial ideal $I$ minimally generated by $\{m_1,\dots,
m_r\}$, we can construct an analogue of the well known
Mayer-Vietoris sequence from topology, in the following way:
\begin{definition}
For each $1\leq s\leq r$ denote $I_s:=\langle m_1,\dots m_s\rangle$, $\tilde{I}_s:=I_{s-1}\cap\langle m_s\rangle=\langle m_{1,s},\dots,m_{s-1,s}\rangle$, where $m_{i,j}$ denotes $\lcm(m_i,m_j)$. Then, for each $s$ we have

\begin{equation}\label{les}
\cdots\longrightarrow H_{i+1}(\KK(I_s))\stackrel{\Delta}{\longrightarrow}H_{i}(\KK(\tilde {I}_s)\longrightarrow H_{i}(\KK(I_{s-1})\oplus \KK(\langle  m_s\rangle  ))\longrightarrow H_{i}(\KK(I_s))\stackrel{\Delta}{\longrightarrow}\cdots
\end{equation}
and since the Koszul differential respects multidegrees, we also
have a multigraded version of the sequence. The set of these
sequences for each $s$ is called the (recursive)
\emph{Mayer-Vietoris sequence} of $I$.
\end{definition}

Using recursively these exact sequences for every $\mu\in\NN^n$ we can compute the Koszul
homology of $I=\langle m_1,\dots,m_r\rangle$.
The involved ideals can be displayed as a tree, the root of which is $I$ and every node $J$ has
as \textit{children} $\tilde J$ on the left and $J'$ on the right (if $J$ is generated by $r$ monomials, $\tilde J$
denotes $\tilde J_r$ and $J'$ denotes $J_{r-1}$).
This is what we call a \textbf{Mayer-Vietoris Tree} of the monomial ideal $I$, and we will denote it $\MVT(I)$. Each
node in a Mayer-Vietoris tree is given a \emph{position}: the root has position $1$ and the left and right children of
the node in position $p$ have respectively, positions $2p$ and $2p+1$. The node of $\MVT(I)$ in position $p$ is denoted $\MVT_p(I)$. We assign also \emph{dimensions} to the nodes: the root has dimension $0$, and the left and right children of
the node of dimension $i$ have respectively, dimensions $i+1$ and $i$.

\begin{remark}
Strictly speaking, the definition of Mayer-Vietoris sequences of monomial ideals is not fully precise, in the sense that the Mayer-Vietoris sequence associated to a given ideal is not uniquely defined, it depends on how the minimal generators are sorted. The choice of the last generator of the ideal $I$ to be the one which defines the {\it Mayer-Vietoris sequence} is just a matter of convenience in notation. Several selection strategies can be applied to select the distinguished generator.
\end{remark}

The properties of Mayer-Vietoris trees allow us to perform Koszul homology computations, see the proofs of the following results in \cite{S08}.

\begin{proposition}
If $H_{i,\mu}(\KK(I))\neq 0$ for some $i$, then $x^{\mu}$ is a generator of some node $J$ in any Mayer-Vietoris tree $\MVT(I)$.
\end{proposition}

Thus, all the multidegrees of Koszul generators (equivalently Betti
numbers) of $I$ appear in $\MVT(I)$. For a sufficient condition, we
need the following notation: among the nodes in $\MVT(I)$ we call
{\it relevant nodes}  those in an even position or in position 1.

\begin{proposition}
If $x^{\mu}$ appears only once as a generator of a relevant node $J$ in $\MVT(I)$ then there exists exactly one generator in $H_*(\KK(I))$ which has multidegree $\mu$.
\end{proposition}

The homological degree to which relevant multidegrees contribute, is just the dimension of the node in which they appear.

Let $I$ be a monomial ideal and $\MVT(I)$ a Mayer-Vietoris tree of
$I$. Let $\mu\in\NN^n$; let $\overline {\beta}_{i,\mu}(I)=1$
if $\mu$ is the multidegree of some non repeated generator in
some relevant node of dimension $i$ in $\MVT(I)$ and $\overline
{\beta}_{i,\mu}(I)=0$ in other case. Let $\widehat
{\beta}_{i,\mu}(I)$ be the number of times $\alpha$ appears as
the multidegree of some generator of dimension $i$ in some relevant
node in $\MVT(I)$. Then for all $\mu\in\NN^n$ we have

$$\overline{\beta}_{i,\mu}(I)\leq\beta_{i,\mu}(I)\leq\widehat {\beta}_{i,\mu}(I).$$

\subsection{Mayer-Vietoris computation of $\Bc_{n-1}$}

If we have a monomial ideal with $r$ minimal generators, we know that $projdim(I)\leq r-1$ (to see this just consider the length of the well known Taylor resolution, for instance), therefore to have $H_{n-1}(I)\neq 0$ we need at least $n$ generators. Using this simple argument, we can prune the Mayer-Vietoris tree of $I$ in the following way: Do not consider all subtrees hanging of nodes with a number of generators less than $n-i$ where $i$ is the dimension of the node. We compute the pruned tree and store the generators of all nodes of dimension $n-1$. This is what we call \emph{pruning by number of generators}.

When we reach a node of $\MVT(I)$ such that the union of the supports of its generators is strictly contained in $\{1,\dots,n\}$ then it is clear that there is no $(n-1)$-st Koszul homology in the subtree hanging from this node, for all maximal corners have full support, therefore, we do not consider this subtree. This is what we call \emph{prunning by number of indeterminates}.

\begin{remark}
Another modification we can apply to the general $\MVT$ algorithm is the following: Consider that we always sort the generators of the nodes in the tree according to the lexicographic term order. Then it is clear that the first generator will have highest exponent in the first indeterminate, we take this as pivot monomial. If the first indeterminate is not $x_1$ we can prune the tree by number of indeterminates. If it is $x_1$ then all generators of the left child of the node will have the same exponent in the first indeterminate, since we take least common multiples for its construction. Therefore, we can consider the node as an ideal in a ring with one fewer indeterminate. Reducing the number of indeterminates, we will have that all our nodes will be considered as ideals in two indeterminates, and here the computation of the $(n-1)$-st Koszul homology is straightforward.
\end{remark}

After performing these prunnings we have a set of multidegrees that are \emph{candidates} to be in $\Bc_{n-1}(I)$. For each of them we check whether they have homology or not based on the caracterization given in examples \ref{ex:algebraic_n-1} and \ref{ex:simplicial_n-1}, i.e.  we just check whether 

$$\frac{x_i\cdot x^\mu}{x_1\cdots x_n}\in I \ \forall i
\quad \mbox{ and } \quad
\frac{x^\mu}{x_1\cdots x_n}\notin I$$

Observe that in the case of squarefree candidates this amounts to check whether they have full support.

\begin{example}
Consider the ideal $I=\langle x^2y^3,y^3zt,yt^2,z^3t^2 \rangle$ in $R=\kb[x,y,z,t]$.
The first level of the lexicographic $\MVT(I)$ is
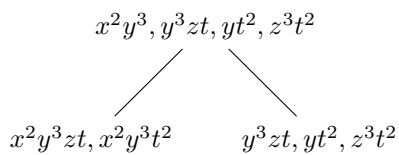
 \begin{center}
\begin{tikzpicture}[scale=1]
\tikzstyle{level 1}=[sibling distance=3cm]
\tikzstyle{level 2}=[sibling distance=2cm]
\node{$x^2y^3,y^3zt,yt^2,z^3t^2$}
child{ node {$x^2y^3zt,x^2y^3t^2$}}
child{ node {$y^3zt,yt^2,z^3t^2$}};
\end{tikzpicture} 
\end{center}
As we can see, the number of generators of $\MVT_2(I)$ is just $2$, so we can prune this branch by number of generators. On the other hand, the $x$ indeterminate is not present in the generators of $\MVT_3(I)$ hence we can prune this branch by number of indeterminates. Therefore, after one iteration of our algorithm we obtain $\Bc_{n-1}(I)=\emptyset$.

\end{example}
\subsection{Comparison with other algorithms}
We present here some tables to show the performance of an implementation of the described algorithm using the {\tt C++} library \cocoalib \cite{cocoalib} which is part of the computer algebra system \cocoa \cite{cocoa}. The Tables show the timings of several different algorithms computing irreducible decompositions of monomial ideals on some benchmark examples. The examples and timings are taken from \cite{R08}, except those corresponding to the Mayer-Vietoris algorithm. Although the tested implementation is not fully optimized to achieve time efficiency, the tables illustrate the fact that the algorithm based on Koszul homology computation has good performance when compared with the algorithms existing in the literature. This is remarkable since the algortihm we present is a simple modification of a more general algorithm, while the others are algorithms specifically designed for the computation of irreducible decomposition. In the Tables $\vert irr(I)\vert$ gives the number of minimal generators of $I$, $\vert irr(I)\vert$ and is the number of irreducible components of $I$. The column \emph{Macaulay2} shows the timings of the computation using the computer algebra system Macaulay2 \cite{M2}, the columns \emph{Monos} correspond to the program Monos by R.A. Milowski \cite{M04} either using the Alexander dual or the Scarf complex approach, and the column \emph{Frobby} corresponds to the \emph{Slice algorithm} by B. Roune \cite{R08}, which shows the best behaviour.

 \begin{table}
 \caption{Random generic ideals in ten indeterminates}
 \label{tab:1}       
 \begin{tabular}{|l|l||l|l|l|l||l|}
 \hline
 $\vert min(I)\vert$ & $\vert irr(I)\vert$ & Macaulay2 & Monos(Alex.)&Monos(Scarf) & Frobby & $\MVT$  \\
 \hline
 40 & 52131 & 226s & 521s & 10 s & 1s & 2.52s \\
 80 & 163162 & OOM & OOT & 54s & 4s & 8.88s\\
 120 & 411997 & OOM & OOT & 198s & 9s & 24.27s\\
 160 & 789687 & OOM &  - & 563s & 19s & 50.78s\\
 200 & 1245139 & OOM & - & OOM & 29s & 86.50s\\
\hline
 \end{tabular}
 \end{table}

\begin{table}
 \caption{Random non-generic ideals in ten indeterminates}
 \label{tab:2}       
 \begin{tabular}{|l|l||l|l|l|l||l|}
 \hline
 $\vert min(I)\vert$ & $\vert irr(I)\vert$ & Macaulay2 & Monos(Alex.)&Monos(Scarf) & Frobby & $\MVT$  \\
 \hline
 100 & 19442 & 43s & 257s & 95s & 1s & 1.96s \\
 150 & 52781 & OOM & 2537s & 3539s & 2s & 7.34s\\
 200 & 79003 & OOM & OOT & 6376s & 3s & 11.53s\\
 400 & 193638 & OOM &  - & OOM & 8s & 36.90s\\
 600 & 318716 & OOM & - & OOM & 16s & 69.88s\\
 800 & 435881 & OOM & - & OOM & 23s & 106.46s\\
 1000 & 571756 & OOM & - & OOM & 32s & 150.04s\\
\hline
 \end{tabular}
 \end{table}


\section{Conclusions and future work}
We have seen in this paper that the Koszul homology of a monomial ideal provides information about the (algebraic) structure of the ideal. In particular we have showed that the $(n-1)$-st Koszul homology of a monomial ideal $I$ provides the irredundant irreducible decomposition of $I$ (or equivalently its Alexander dual, set of maximal standard monomials,...) and a Stanley decomposition of $R/I$. These facts allow us to produce algorithm to make computations on monomial ideals based on the computation of its Koszul homology. This has been illustrated by the computation of the irreducible decomposition of a monomial ideal using a specialization of the \emph{Mayer-Vietoris tree algorithm}, showing good performance in practice.

Future work includes the exploration of further relations between algebraic and homological descriptions and properties of monomial ideals, based on their combinatorial nature. Also the application of this homological approach to other algebraic combinatoric objects is an open line of research already followed by different authors, see as good examples \cite{MS04} and \cite{OW07}. Finally, the improvement of the algorithms presented in this paper and the implementation of new algorithms to perform computations on monomial ideals is another direction for future work already in progress by the authors.

\noindent{\bf Aknowledgements:}
The authors want to thank Werner M. Seiler for his advise and help, and John Abbott for his help in the implementation issues of the algorithms in section \ref{sec:algorithms}. Also the second author wants to thank the Dipartimento di Matematica of the University of Genova (Italy) for their hospitality during part of this work. The second author was partially supported by project Fomenta 03/07 of the \emph{Comunidad Autonoma de La Rioja} (Spain).

\bibliographystyle{plain} 
\bibliography{n-1_Koszul_homology}   


\end{document}